\documentclass[11pt,reqno,oneside]{amsart}
\usepackage{verbatim,amsmath,amssymb,cite,epsfig,xspace,color}

\numberwithin{equation}{section}

\theoremstyle{plain}
\newtheorem{theorem}{Theorem}[section]

\theoremstyle{definition}

\theoremstyle{remark}
\newtheorem{remark}{Remark}[section]

\newcommand{\ignore}[1]{}

\newcommand{\yb}{\mathbf y}

\newcommand{\E}{{\mathcal E}}

\newcommand{\la}{\langle}
\newcommand{\ra}{\rangle}

\newcommand{\half}{\frac{1}{2}}

\newcommand{\eps}{\epsilon}

\begin{document}

\title
[Lattice Stability for Atomistic Chains Modeled by the EAM] {Lattice Stability
for Atomistic Chains Modeled by Local Approximations of the Embedded Atom Method
}
\author{Xingjie Helen Li and Mitchell Luskin}


\thanks{
This work was supported in part by DMS-0757355,
 DMS-0811039,  the Institute for Mathematics and
Its Applications, and
 the University of Minnesota Supercomputing Institute.
  This work was also supported by the Department of Energy under Award Number
DE-SC0002085.
}



\date{\today}

\begin{abstract}
The accurate approximation of critical strains for lattice instability
is a key criterion for predictive computational modeling of materials.  In this
paper, we present a comparison of the lattice stability for atomistic chains
modeled by the embedded atom method (EAM) with their approximation
by local Cauchy-Born models.  We find that both the volume-based local model
and the reconstruction-based local model can give O$(1)$ errors for the
critical strain since the embedding energy density is generally strictly convex.
The critical strain predicted by the volume-based model is
always larger than that predicted by the atomistic model, but the critical
strain for reconstruction-based models can be either larger or smaller than
that predicted by the atomistic model.
\end{abstract}

\maketitle{\thispagestyle{empty} \maketitle
\section{Introduction}

Predictive multiscale computational methods must be accurate
near lattice instabilities that characterize the formation and movement
of cracks, dislocations, and grain boundaries.
In this paper, we present analytic results comparing the
lattice instabilities predicted by an atomistic chain modeled
by the embedded atom method (EAM) with the lattice instabilities
predicted by local approximations of the atomistic model.

Since it is not possible to compute large enough
fully atomistic systems to accurately approximate the interaction
of local defects with long-range elastic fields, atomistic-to-continuum
coupling methods have been proposed
~\cite{curt03,LinP:2006a,Miller:2003a,Shimokawa:2004,E:2004,Miller:2008,Legoll:2005,shapeev,bqce11}.
For crystalline solids, the
continuum region is generally computed by coarse-graining a
{\em local} approximation of the {\em nonlocal} atomistic model.
An atom in the nonlocal atomistic region interacts with all
of its neighbors within a cutoff radius.  In the continuum region,
the Cauchy-Born rule is used to derive a local model that
approximates the interactions of
atoms beyond their nearest neighbors by modified interactions
of their nearest neighbors.

To verify that an atomistic-to-continuum coupling method
accurately reproduces the lattice stability of the fully
atomistic model, it is necessary to first verify that
the local (continuum) model itself reproduces the lattice
stability of the fully atomistic model.  Even if the
local model reproduces the lattice stability of the
fully atomistic model, the atomistic-to-continuum coupling method
may not reproduce the lattice stability of the
fully atomistic model because of the error introduced by the
coupling~\cite{doblusort:qce.stab,luskin.durham}.

It has been proven in ~\cite{doblusort:qce.stab} that the Cauchy-Born
local model reproduces the lattice stability for an atomistic
chain modeled by Lennard-Jones type pair interaction
(we note that the volume-based and reconstruction-based local models
are equal for a pair potential interaction).  And it has been proven
for multidimensional lattices that the set of stable uniform strains
for the atomistic model is a subset of the set of uniform strains
for the Cauchy-Born volume-based local model~\cite{HudsOrt:a,E:2007a}, but the
equality of these sets has not been demonstrated analytically.
As a matter of fact, numerical experiments in~\cite{HudsOrt:a}
suggest that the the inclusion is strict in some cases.

In this paper, we prove for an atomistic chain that not only
are the sets of stable uniform strains different for the atomistic model
and the local models for a many-body potential, but the set of stable
uniform strains can be different for {\em volume-based} and {\em reconstruction-based}
local models.
We will focus our analysis on the embedded atom method~\cite{Foiles,Mishin,Johnson:1989},
which is an empirical many-body
potential that is widely used to model FCC metals such as copper
and aluminum.
We identify the critical assumptions for the pair potential, electron density
function, and embedding function to study the lattice stability of the atomistic and the different
local models.  We find that both the volume-based local model
and the reconstruction-based local model can give O$(1)$ errors for the
critical strain since the embedding energy density is generally strictly convex.

In Section~\ref{notation}, we present the notation used in this
paper. We define the displacement space $\mathcal{U}$ and the
deformation space $\mathcal{Y}_{F}$. We then introduce the norms we
will use to estimate the modeling error and the displacement
gradient error.
In Section~\ref{QCapproximations}, we briefly review the formulae of
the fully atomistic EAM model and the volume-based and the
reconstruction-based local quasicontinuum (QCL) model, respectively.

In Section~\ref{sharp}, we give precise stability estimates for
the fully atomistic model, the volume-based and the
reconstruction-based local models for a uniformly strained
chain. We then compare the stability conditions of each model
under different assumptions.  We summarize our results and
discuss extensions to multidimensional issues in the Conclusion.


\section{Notation}
\label{notation}
In this section, we
present the notation used in this paper. We define the scaled reference
lattice
\[
\eps \mathbb{Z}:= \{\eps\ell : \ell\in\mathbb{Z}\},
\]
 where $\eps >0$
scales the reference atomic spacing and $\mathbb{Z}$ is the set of
integers.  We then deform the reference lattice $\eps \mathbb{Z}$
uniformly into the lattice
\[
F\eps \mathbb{Z}:= \{F\eps\ell : \ell\in\mathbb{Z}\}
\]
where $F >0$ is the macroscopic deformation gradient, and we define
the corresponding deformation $\mathbf{y}_{F}$ by
\[
(\mathbf{y}_{F})_{\ell}:=F\eps \ell  \quad \text{for } -\infty
<\ell<\infty.
\]
For simplicity, we consider the space $\mathcal{U}$ of $2N$-periodic
zero mean displacements $\mathbf{u}=(u_{\ell})_{\ell \in
\mathbb{Z}}$ from $\mathbf{y}_{F}$ given by
\[
\mathcal{U}:=\bigg\{\mathbf{u} : u_{\ell+2N}=u_{\ell}
\text{ for }\ell\in \mathbb{Z},\,
\text{and}\sum_{\ell=-N+1}^{N}u_{\ell}=0\bigg\},
\]
and we thus admit deformations $\mathbf{y}$ from the space
\[
\mathcal{Y}_{F}:=\{\mathbf{y}:
\mathbf{y}=\mathbf{y}_{F}+\mathbf{u}\text{ for some }\mathbf{u}\in
\mathcal{U}\}.
\]
We set $\eps=1/N$ throughout so that the reference length of the
periodic domain is fixed.

 We
define the discrete differentiation operator, $D\mathbf{u}$, on
periodic displacements by
\[
(D\mathbf{u})_{\ell}:=\frac{u_{\ell}-u_{\ell-1}}{\epsilon}, \quad
-\infty<\ell<\infty.
\]
We note that $\left(D\mathbf{u}\right)_{\ell}$ is also $2N$-periodic
in $\ell$ and satisfies the zero mean condition. We will denote
$\left(D\mathbf{u}\right)_{\ell}$ by $Du_{\ell}$.
We then define
\begin{align*}
\left(D^{(2)}\mathbf{u}\right)_{\ell}:=\frac{Du_{\ell}-Du_{\ell-1}}{\epsilon},
\qquad -\infty<\ell<\infty,
\end{align*}
and we define $\left(D^{(3)}\mathbf{u}\right)_{\ell}$ and
$\left(D^{(4)}\mathbf{u}\right)_{\ell}$ in a similar way. To make
the formulas concise and more readable, we sometimes denote $Du_{\ell}$ by $u'_{\ell}$,
$D^{(2)}u_{\ell}$ by $u''_{\ell}$, etc., when there is no confusion
in the expressions.

For a displacement $\mathbf{u}\in \mathcal{U}$ and its discrete derivatives, we define the
discrete $\ell_{\epsilon}^{2}$ norms by
\begin{align*}
\|\mathbf{u}\|_{\ell_{\epsilon}^{2}}&:= \left( \epsilon
\sum_{\ell=-N+1}^{N}|u_{\ell}|^{2}\right)^{1/2},\qquad
\|\mathbf{u}'\|_{\ell_{\epsilon}^{2}}:= \left( \epsilon
\sum_{\ell=-N+1}^{N}|u_{\ell}'|^{2}\right)^{1/2},\text{ etc.}
\end{align*}
Finally, for smooth real-valued functions $\E(\yb)$ defined for
$\yb\in\mathcal{Y}_{F},$ we define the first and second derivatives (variations) by
\begin{equation*}
\begin{split}
\la\delta\mathcal{E}(\mathbf{y}),\mathbf{w}\ra&:=\sum_{\ell=-N+1}^{N}
 \frac{\partial \mathcal{E}}{\partial y_{\ell}}(\mathbf{y})w_{\ell}\qquad
\text{for all }\mathbf{w}\in \mathcal{U},\\
\la\delta^2\mathcal{E}(\mathbf{y})\mathbf{v},\mathbf{w}\ra&:=\sum_{\ell,\,m=-N+1}^{N}
 \frac{\partial^2 \mathcal{E}}{\partial y_{\ell}\partial y_{m}}(\mathbf{y})v_{\ell}
 w_{m}\qquad
\text{for all }\mathbf{v},\,\mathbf{w}\in \mathcal{U}.
\end{split}
\end{equation*}

\section{The Embedded Atom Model and Its Local Approximations.}\label{QCapproximations}
In this section, we will give a short description for the next-nearest neighbor atomistic EAM model
and its approximations .
\subsection{The Next-Nearest Neighbor Atomistic EAM Model}
Given deformations $\mathbf{y}\in \mathcal{Y}_{F}$,
the total energy per period of the next-nearest neighbor atomistic EAM model is
\begin{equation}\label{AtomTotalEnergy}
 \mathcal{E}_{tot}^{a}(\mathbf{y}):=\mathcal{E}^{a}(\mathbf{y})+\mathcal{F}(\mathbf{y}),
\end{equation}
where $\mathcal{E}^{a}(\mathbf{y})$ is the
total atomistic energy and $\mathcal{F}(\mathbf{y})$
is the total external potential energy.
The total atomistic energy $\mathcal{E}^{a}(\mathbf{y})$ is the sum of the {\it embedding energy}, $\hat{\mathcal{E}}^{a}(\mathbf{y}),$
and the {\it pair potential energy}, $\tilde{\mathcal{E}}^{a}(\mathbf{y})$. The energy expression is
\begin{equation}\label{AtomModel}
\mathcal{E}^{a}(\mathbf{y}):=\hat{\mathcal{E}}^{a}(\mathbf{y})+\tilde{\mathcal{E}}^{a}(\mathbf{y})
=\epsilon\sum_{\ell=-N+1}^{N}\left(\hat{\mathcal{E}}^{a}_{\ell}(\mathbf{y})+\tilde{\mathcal{E}}^{a}_{\ell}(\mathbf{y})\right).
\end{equation}
The embedding energy per atom (per atomistic reference spacing $\eps$) is defined as
$
\hat{\mathcal{E}}^{a}_{\ell}(\mathbf{y}):=
                            G\left(\bar\rho^a_\ell(\mathbf{y})\right),
$
where $G(\bar\rho)$ is the embedding energy function and
$\bar\rho^a_\ell(\mathbf{y})$ is the total electron density
at
atom $\ell$:
\[
\bar\rho^a_\ell(\mathbf{y}):=
\rho(y'_{\ell})+\rho(y'_{\ell}+y'_{\ell-1})
                            +\rho(y'_{\ell+1})+\rho(y'_{\ell+1}+y'_{\ell+2}).
\]
The function $\rho(r/\eps)$ is the electron density contributed by an atom at distance $r.$

The pair potential energy per atom (per atomistic reference spacing $\eps$) is
\[
 \tilde{\mathcal{E}}^{a}_{\ell}(\mathbf{y}):=
 \half \left[\phi(y'_{\ell})+\phi(y'_{\ell}+y'_{\ell-1}) +\phi(y'_{\ell+1})+\phi(y'_{\ell+1}+y'_{\ell+2})\right],
\]
where $\phi(r/\eps)$ is
the pair potential interaction energy~\cite{Foiles}.
Our formulation
allows general nonlinear external potential energies $\mathcal{F}(\mathbf{y})$ defined
for $\mathbf{y}\in \mathcal{Y}_{F}$, but for simplicity, we only consider the
total external potential
energy for $2N$-periodic dead loads $\mathbf{f}$
\[
 \mathcal{F}(\mathbf{y}):=-\sum_{\ell=-N+1}^{N}\epsilon f_{\ell}y_{\ell}.
\]
The equilibrium solution $\mathbf{y}^{a}$ of the EAM-atomistic
model \eqref{AtomTotalEnergy} then satisfies
\begin{align}\label{AtomEqulibriumEq2}
- \la\delta \mathcal{E}^{a}(\mathbf{y}^a),\mathbf{w}\ra=
-\la\delta\hat{\mathcal{E}}^{a}(\mathbf{y}^a),\mathbf{w}\ra
-\la\delta\tilde{\mathcal{E}}^{a}(\mathbf{y}^a),\mathbf{w}\ra
=\la\delta\mathcal{F}(\mathbf{y}^{a}),\mathbf{w}\ra \qquad \text{for all }
\mathbf{w}\in \mathcal{U}.
 \end{align}
Here the negative of the embedding force of \eqref{AtomEqulibriumEq2} is
\begin{align*}
\la\delta\hat{\mathcal{E}}^{a}(\mathbf{y}^a),\mathbf{w}\ra&=\epsilon
\sum_{\ell=-N+1}^{N}G'\Big(\bar\rho^a_\ell(\mathbf{y}^a) \Big)\cdot\Big[
\rho'(Dy^{a}_{\ell})w'_{\ell}+\rho'(Dy^{a}_{\ell}+Dy^{a}_{\ell-1})(w'_{\ell}+w'_{\ell-1})\Big.\\
&\qquad\qquad\qquad
\Big.+\rho'(Dy^{a}_{\ell+1})w'_{\ell+1}+\rho'(Dy^{a}_{\ell+1}+Dy^{a}_{\ell+2})(w'_{\ell+1}+w'_{\ell+2})
\Big],
\end{align*}
the negative of the pair potential force of \eqref{AtomEqulibriumEq2} is given by
\begin{align*}
\la\delta\tilde{\mathcal{E}}^{a}(\mathbf{y}^a),\mathbf{w}\ra&=
\epsilon\sum_{\ell=-N+1}^{N}\half
\Big[
\phi'(Dy^{a}_{\ell})w'_{\ell}+\phi'(Dy^{a}_{\ell}+Dy^{a}_{\ell-1})(w'_{\ell}+w'_{\ell-1})
\\
&\qquad\qquad\qquad\qquad\qquad
+\phi'(Dy^{a}_{\ell+1})w'_{\ell+1}+\phi'(Dy^{a}_{\ell+1}+Dy^{a}_{\ell+2})(w'_{\ell+1}+w'_{\ell+2})
\Big],
\end{align*}
and the negative of the external force is formulated as
\[
 \la\delta\mathcal{F}(\mathbf{y}),\mathbf{w}\ra=\sum_{\ell=-N+1}^{N}
 \frac{\partial \mathcal{F}}{\partial y_{\ell}}(\mathbf{y})w_{\ell}
 =-\sum_{\ell=-N+1}^{N}\epsilon f_{\ell}w_{\ell}.
\]
\subsection{The Local EAM Approximations.}
In this subsection, we will briefly review the idea of the two different local
 approximations, the volume-based and the
reconstruction-based, and give their expressions respectively.
\subsubsection{The Volume-Based Local EAM Approximation.}
The idea of the volume-based local approximation based on the
Cauchy-Born rule was first proposed in \cite{Ortiz:1995a,Shenoy:1999a,Miller:2003a}.
We denote this energy by $\mathcal{E}^{c,v}(\mathbf{y}),$
and we can formulate the local energy associated with each
atom as
\begin{align*}
 \mathcal{E}^{c,v}_{\ell}(\mathbf{y}):=
 \hat{\mathcal{E}}^{c,v}_{\ell}(\mathbf{y})+\tilde{\mathcal{E}}^{c,v}_{\ell}(\mathbf{y})
=&
\half G\left(\bar\rho^{c,v}_\ell(\mathbf{y}))\right)+\half
G\left(\bar\rho^{c,v}_{\ell+1}(\mathbf{y})\right) \\
&\quad
+\half \left[\phi(y'_{\ell})+\phi(2y'_{\ell})
+\phi(y'_{\ell+1})+\phi(2y'_{\ell+1})\right],
\end{align*}
where the total local electron density at atom $\ell$ is
\[
\bar\rho^{c,v}_\ell(\mathbf{y}):=2\rho(y'_{\ell})
 +2\rho(2y'_{\ell}).
\]
Then the total volume-based local  energy is
\begin{align}\label{QCLenergy1}
\mathcal{E}^{c,v}_{tot}(\mathbf{y}):=
\mathcal{E}^{c,v}(\mathbf{y})+\mathcal{F}(\mathbf{y})=\epsilon\sum_{\ell=-N+1}^{N}\mathcal{E}^{c,v}_{\ell}(\mathbf{y})
-\epsilon\sum_{\ell=-N+1}^{N}f_{\ell}y_{\ell}.
\end{align}
The equilibrium solution $\mathbf{y}^{c,v}$ then satisfies
\begin{align}\label{QCLEquilibrium1}
-\la\delta\mathcal{E}^{c,v}(\mathbf{y}^{c,v}),\mathbf{w}\ra
=-\la\delta\hat{\mathcal{E}}^{c,v}(\mathbf{y}^{c,v}),\mathbf{w}\ra
-\la\delta\tilde{\mathcal{E}}^{c,v}(\mathbf{y}^{c,v}),\mathbf{w}\ra
=\la\delta\mathcal{F}(\mathbf{y}^{c,v}),\mathbf{w}\ra\quad \text{for all}\, \mathbf{w}\in \mathcal{U}.
\end{align}
The negative of the embedding force of \eqref{QCLEquilibrium1} is
\begin{align*}
\la\delta\hat{\mathcal{E}}^{c,v}(\mathbf{y}^{c,v}),\mathbf{w}\ra
=&\epsilon\sum_{\ell=-N+1}^{N}
\left\{G'\left(\bar{\rho}_{\ell}^{c,v}(\mathbf{y}^{c,v})\right)
\cdot\left[\rho'(Dy^{c,v}_{\ell})+2\rho'(2Dy^{c,v}_{\ell})\right]w'_{\ell}\right.\\
&\qquad+\left.G'\left(\bar{\rho}_{\ell+1}^{c,v}(\mathbf{y}^{c,v})\right)
\cdot\left[\rho'(Dy^{c,v}_{\ell+1})+2\rho'(2Dy^{c,v}_{\ell+1})\right]w'_{\ell+1}\right\},
\end{align*}
and the negative of the pair potential force of \eqref{QCLEquilibrium1} is
given by
\begin{align*}
\la\delta\tilde{\mathcal{E}}^{c,v}(\mathbf{y}^{c,v}),\mathbf{w}\ra
=&\epsilon\sum_{\ell=-N+1}^{N}
\half\left\{\left[\phi'(Dy^{c,v}_{\ell})+2\phi'(2Dy^{c,v}_{\ell})\right]w'_{\ell}
 + \left[\phi'(Dy^{c,v}_{\ell+1})+2\phi'(2Dy^{c,v}_{\ell+1})\right]w'_{\ell+1}\right\}.
\end{align*}
\subsubsection{The Reconstruction-Based Local EAM Approximation.}
Using the Cauchy-Born approximation, one can also reconstruct
 the position of each atom~\cite{E:2004} and compute the energy
 $\mathcal{E}^{c,r}(\mathbf{y})$
 by the approximation
\begin{align*}
 \mathcal{E}^{c,r}_{\ell}(\mathbf{y})=
 \hat{\mathcal{E}}^{c,r}_{\ell}(\mathbf{y})+\tilde{\mathcal{E}}^{c,r}_{\ell}(\mathbf{y})
=&  G\left(\bar\rho^{c,r}_{\ell}(\mathbf{y})\right) +\half \left[\phi(y'_{\ell})+\phi(2y'_{\ell})
+\phi(y'_{\ell+1})+\phi(2y'_{\ell+1})\right],
\end{align*}
where the reconstruction-based local electron density at atom $\ell$ is
\[
\bar\rho^{c,r}_\ell(\mathbf{y}):=\rho(y'_{\ell})+\rho(2y'_{\ell})+
\rho(y'_{\ell+1})+\rho(2y'_{\ell+1}).
\]
Thus, the total energy of the reconstruction-based local
model is
\begin{align}\label{QCLenergy2}
\begin{split}
\mathcal{E}_{tot}^{c,r}(\mathbf{y}):=&\hat{\mathcal{E}}^{c,r}(\mathbf{y})
+\tilde{\mathcal{E}}^{c,r}(\mathbf{y})+\mathcal{F}(\mathbf{y})\\
=&\epsilon\sum_{\ell-N+1}^{N}\Big\{G\left[\rho(y'_{\ell})+\rho(2y'_{\ell})+
\rho(y'_{\ell+1})+\rho(2y'_{\ell+1})\right]\Big. \\
&\qquad  \qquad\Big.+\half \left[\phi(y'_{\ell})+\phi(2y'_{\ell})
+\phi(y'_{\ell+1})+\phi(2y'_{\ell+1})\right]\Big\}-\epsilon\sum_{\ell=-N+1}^{N}f_{\ell}y_{\ell}.
\end{split}
\end{align}
The volume-based and reconstruction-based local energies have the same pair potential energy,
but their approximations for the embedding energy are quite different.

We compute the equilibrium solution of the reconstruction-based
local model \eqref{QCLenergy2} from
\begin{align}\label{QCLEquilibrium2}
-\la\delta\mathcal{E}^{c,r}(\mathbf{y}^{c,r}),\mathbf{w}\ra
=-\la\delta\hat{\mathcal{E}}^{c,r}(\mathbf{y}^{c,r}),\mathbf{w}\ra
-\la\delta\tilde{\mathcal{E}}^{c,r}(\mathbf{y}^{c,r}),\mathbf{w}\ra
=\la\delta\mathcal{F}(\mathbf{y}^{c,r}),\mathbf{w}\ra\quad \text{for all}\ \mathbf{w}\in \mathcal{U}.
\end{align}
Here the negative of the embedding force of \eqref{QCLEquilibrium2} is
\begin{align*}
\la\delta\hat{\mathcal{E}}^{c,r}(\mathbf{y}^{c,r}),\mathbf{w}\ra
=\epsilon&\sum_{\ell=-N+1}^{N}G'\left(\bar{\rho}^{c,r}_{\ell}(\mathbf{y}^{c,r})\right)
\cdot\left[\left(\rho'(Dy^{c,r}_{\ell})+2\rho'(2Dy^{c,r}_{\ell})\right)w'_{\ell}\right.\\
&\qquad\qquad\qquad\qquad\left.
+\left(\rho'(Dy^{c,r}_{\ell+1})+2\rho'(2Dy^{c,r}_{\ell+1})\right)w'_{\ell+1}\right],
\end{align*}
and the negative of the pair potential force of \eqref{QCLEquilibrium2} is
\begin{align*}
\la\delta\tilde{\mathcal{E}}^{c,r}(\mathbf{y}^{c,r}),\mathbf{w}\ra
=&\epsilon\sum_{\ell=-N+1}^{N}
\half\left\{\left[\phi'(Dy^{c,r}_{\ell})+2\phi'(2Dy^{c,r}_{\ell})\right]w'_{\ell}
 + \left[\phi'(Dy^{c,r}_{\ell+1})+2\phi'(2Dy^{c,r}_{\ell+1})\right]w'_{\ell+1}\right\}.
\end{align*}
The pair potential energy of both local approximations are exactly the same, but
the embedding parts are quite different, which leads to different critical
strains for lattice instability.
We will analyze the lattice stability for all of the models in the next section.
\section{Sharp Stability Analysis of The Atomistic and Local EAM Models.}\label{sharp}
In this section, we analyze and compare the
conditions for lattice stability of the atomistic model and the two local
approximations for the next-nearest neighbor case. We will use
techniques similar to those presented in \cite{doblusort:qce.stab}
for the atomistic and  quasicontinuum methods with pair potential
interaction.
\subsection{Stability of the Atomistic EAM Model.}
We first consider the fully atomistic model.
The uniform deformation $\mathbf{y}_{F}$ is an equilibrium of the
atomistic model~\eqref{AtomModel} without external force. We call
$\mathbf{y}_{F}$ stable in the atomistic model if and only if
$\la\delta^{2}\mathcal{E}^{a}(\mathbf{y}_{F})$ is positive definite,
that is,
\begin{align}\label{AtomStabEq0}
 \la\delta^{2}\mathcal{E}^{a}(\mathbf{y}_{F})\mathbf{u},\mathbf{u}\ra=
\la\delta^{2}\hat{\mathcal{E}}^{a}(\mathbf{y}_{F})\mathbf{u},\mathbf{u}\ra
+\la\delta^{2}\tilde{\mathcal{E}}^{a}(\mathbf{y}_{F})\mathbf{u},\mathbf{u}\ra>0
\quad \text{for all } \mathbf{u}\in \mathcal{U}\setminus \{\mathbf{0}\}.
\end{align}
We computed
$\la\delta^{2}\tilde{\mathcal{E}}^{a}(\mathbf{y}_{F})\mathbf{u},\mathbf{u}\ra$
 in \cite{doblusort:qce.stab} to obtain
\begin{align}\label{PairAtomStabEq1}
 \begin{split}
\la\delta^{2}\tilde{\mathcal{E}}^{a}(\mathbf{y}_{F})\mathbf{u},\mathbf{u}\ra
=\tilde{A}_{F}\|D\mathbf{u}\|_{\ell_{\epsilon}^2}^2
-\epsilon^2\phi''_{2F}\|D^{(2)}\mathbf{u}\|_{\ell_{\epsilon}^2}^2,
 \end{split}
\end{align}
where
\begin{equation}\label{PairAtomCond1}
\tilde{A}_{F}:=\phi''_{F}+4\phi''_{2F}\quad\text{for}\quad
\phi''_{F}:=\phi''(F)\text{ and } \phi''_{2F}:=\phi''(2F)
\end{equation}
is the {\it continuum elastic modulus for the pair interaction potential}.
 Thus,
we focus on
$\la\delta^{2}\hat{\mathcal{E}}^{a}(\mathbf{y}_{F})\mathbf{u},\mathbf{u}\ra$,
which can be formulated as
\begin{align}\label{AtomStabEq1}
\begin{split}
 \la\delta^{2}\hat{\mathcal{E}}^{a}(\mathbf{y}_{F})\mathbf{u},\mathbf{u}\ra
&=\epsilon \sum_{\ell=-N+1}^{N}\Bigg\{G''_{F}\,
                          \left[\rho'_{F}(u'_{\ell}+u'_{\ell+1})
                         +\rho'_{2F}(u'_{\ell-1}+u'_{\ell}+u'_{\ell+1}+u'_{\ell+2})\right]^2\\
&\qquad\qquad\qquad\quad\left. +G'_{F}\left[\rho''_{F}(u'_{\ell})^2+\rho''_{2F}(u'_{\ell}+u'_{\ell-1})^2
              +\rho''_{F}(u'_{\ell+1})^2\right.\right.\\
&\qquad\qquad\qquad\qquad\qquad\quad\left.+\rho''_{2F}(u'_{\ell+1}+u'_{\ell+2})^2\right]\Bigg\},
\end{split}
\end{align}
where we use the simplified notation
\begin{gather*}
\rho'_{F}:=\rho'(F),\quad \rho''_{F}:=\rho''(F),
\quad \rho'_{2F}:=\rho(2F),\quad \rho''_{2F}:=\rho''(2F),\\
 G'_{F}:=G'(\bar\rho^a_\ell(\mathbf{y}_F))=G'(\bar\rho^{c,v}_\ell(\mathbf{y}_F))
 =G'(\bar\rho^{c,r}_\ell(\mathbf{y}_F)),\\
 G''_{F}:=G''(\bar\rho^a_\ell(\mathbf{y}_F))
 =G''(\bar\rho^{c,v}_\ell(\mathbf{y}_F))
 =G''(\bar\rho^{c,r}_\ell(\mathbf{y}_F)).
\end{gather*}

We define the {\it continuum elastic modulus for the embedding energy} to be
\begin{equation}\label{EAMCond1}
\hat{A}_{F}:=4G''_{F}\left(\rho'_{F}+2\rho'_{2F}\right)^2
+2G'_{F}\left(\rho''_{F}+4\rho''_{2F}\right),
\end{equation}
and we define
\begin{gather}\label{AtomCoeff}
\begin{split}
A_{F}:=\hat{A}_{F}+\tilde{A}_{F},\quad
B_{F}:= -\left[\phi''_{2F}+G''_{F}\Big((\rho'_{F})^2+20(\rho'_{2F})^2+12\rho'_{F}\rho'_{2F}\Big)+G'_{F}\left(2\rho''_{2F}\right)\right],\\
C_{F}:=G''_{F}\left(8(\rho'_{2F})^2+2\rho'_{F}\rho'_{2F}\right),\quad\text{and}\quad
 D_{F}:=-G''_{F}\left(\rho'_{2F}\right)^2.
\end{split}
\end{gather}
Then \eqref{AtomStabEq0} becomes
\begin{equation}\label{Four}
\begin{split}
\la\delta^{2}\mathcal{E}^{a}(\mathbf{y}_{F})\mathbf{u},\mathbf{u}\ra=&
A_F\|D\mathbf{u}\|_{\ell_{\epsilon}^{2}}^{2}
+\epsilon^2B_F \|D^{(2)}\mathbf{u}\|_{\ell_{\epsilon}^{2}}^{2}
+\epsilon^4
C_F
\|D^{(3)}\mathbf{u}\|_{\ell_{\epsilon}^2}^{2}
+\epsilon^6 D_F
\|D^{(4)}\mathbf{u}\|_{\ell_{\epsilon}^2}^{2},
\end{split}
\end{equation}
where the detailed calculation can be found in the paper \cite{LuskinXingjieEAM:2010}.

We will analyze the stability of $\la\delta^{2}\mathcal{E}^{a}(\mathbf{y}_{F})\mathbf{u},\mathbf{u}\ra$
by using the Fourier representation~\cite{HudsOrt:a}
\begin{equation}\label{fourier}
Du_{\ell}=\sum_{\substack {
k=-N+1\\ k\neq 0}}
^{N}\frac{c_{k}}{\sqrt{2}}\cdot\exp\left(i\,k\frac{\ell}{N}\pi\right).
\end{equation}
We exclude $k=0$ since $D\mathbf{u}$ must satisfy the mean zero
condition $\sum_{\ell=-N+1}^{N}Du_{\ell}=0$.

It then follows from the discrete orthogonality of the Fourier basis
that
\begin{equation}\label{AtomStabFourierEq1}
\begin{split}
\la\delta^{2}\mathcal{E}^{a}(\mathbf{y}_{F})\mathbf{u},\mathbf{u}\ra &=
\sum_{\substack {
k=-N+1\\ k\neq 0}}
^{N} |c_{k}|^2\cdot \Bigg\{ A_F +B_F
\left[4\sin^{2}\left(\frac{k\pi}{2N}\right)\right]\\
&\qquad\qquad\qquad +C_F\left[4\sin^2\left(\frac{k\pi}{2N}\right)\right]^2
+D_F\left[4\sin^2\left(\frac{k\pi}{2N}\right)\right]^3\Bigg\}.
\end{split}
\end{equation}
We see from \eqref{AtomStabFourierEq1} that the eigenvalues
$\lambda^a_k$ for $k=1,\dots,N$ of
$\la\delta^{2}\mathcal{E}^{a}(\mathbf{y}_{F})\mathbf{u},\mathbf{u}\ra$
with respect to the $\|D\mathbf{u}\|_{\ell_{\epsilon}^{2}}$ norm are
given by
\[
\lambda^a_k=\lambda^a_{F}(s_k)\quad\text{for }\quad s_k=4\sin^2\left(\frac{k\pi}{2N}\right)
\]
where
\[
\lambda^a_{F}(s):=A_{F}+B_{F}s+C_{F}s^2+D_{F}s^3.
\]

The energy and electron densities figures in \cite{Foiles} and \cite{Mishin}
satisfy the following conditions which we shall assume in our analysis
\begin{align}\label{EAMfunAssumption1}
\begin{split}
\phi''_{F}>0,\,\phi''_{2F}<0;\quad
 \rho'_{F}\le 0, \, \rho'_{2F} \le 0;\quad
 \rho''_{F}\ge 0,\,\rho''_{2F} \ge 0;\quad\text{and}\quad
G''_{F}\ge 0.
\end{split}
\end{align}
We can derive from the assumption \eqref{EAMfunAssumption1} that
\begin{align}\label{condition1}
C_{F}> 0,\quad
D_{F}< 0,\quad\text{and}\quad
8|D_{F}|\le C_{F}.\quad
\end{align}
Since \eqref{condition1} implies that $|D_Fs|\le 4|D_F|\le C_F/2,$ for $0\le s\le 4,$ we have that
\begin{equation}\label{lambda}
{\lambda^a_{F}}'(s)=B_F+2C_Fs+3D_Fs^2\ge B_{F}+\frac {C_{F}}2 s\quad\text{for all}\quad 0\le s\le 4.
\end{equation}

We note from \eqref{lambda} that the condition $B_F\ge 0$ or equivalently
\begin{align}\label{EAMfunAssumption2}
 \begin{split}
&\phi''_{2F}+G''_{F}\left[\left(\rho'_{F}\right)^2
+20\left(\rho'_{2F}\right)^2+12\rho'_{F}\rho'_{2F}\right]
+G'_{F}\,2\rho''_{2F}=-B_F\le 0,
 \end{split}
\end{align}
implies that $\lambda^{a}_{F}(s)$ is increasing for
$0\le s \le 4.$  We thus conclude that
if $B_{F}\ge 0,$ then
\begin{equation}\label{eq:sharp}
\la\delta^{2}\mathcal{E}^{a}(\mathbf{y}_{F})\mathbf{u},\mathbf{u}\ra
\ge {\lambda^a_{F}}(s_1) \|D\mathbf{u}\|_{\ell_{\epsilon}^{2}}^2\ge \left(\hat{A}_{F}+\tilde{A}_{F}\right)
\|D\mathbf{u}\|_{\ell_{\epsilon}^{2}}^2\quad\text{for all }\mathbf{u}\in \mathcal{U}.
\end{equation}
This result is summarized in the following theorem:
\begin{theorem}\label{AtomStabThm}
Suppose that the hypotheses \eqref{EAMfunAssumption1} and $B_{F}\ge 0$ hold.
Then the uniform deformation $\mathbf{y}_{F}$ is stable for the
atomistic model if and only if
\begin{align*}
\lambda^a_{F}(s_1)&= A_F
+B_F
\left[4\sin^{2}\left(\frac{\pi}{2N}\right)\right]
+C_F\left[4\sin^2\left(\frac{\pi}{2N}\right)\right]^2
+D_F\left[4\sin^2\left(\frac{\pi}{2N}\right)\right]^3\\
  &=\hat{A}_{F}+\tilde{A}_{F}-4\sin^2\left(\frac{\pi}{2N}\right)\left\{\phi''_{2F}+G''_{F}\left[\left(\rho'_{F}\right)^2
+20\left(\rho'_{2F}\right)^2+12\rho'_{F}\rho'_{2F}\right]
+G'_{F}\,2\rho''_{2F}\right\}\\
&\qquad+4^2\sin^4\left(\frac{\pi}{2N}\right)G''_{F}
\left[\eta\left(\rho'_{2F}\right)^2+2\rho'_{F}\rho'_{2F}\right]
-4^3\sin^6\left(\frac{\pi}{2N}\right)G''_{F}\left(\rho'_{2F}\right)^2>0.
\end{align*}
\end{theorem}
We note that the differences between $s_{k}$ and $s_{k-1}$ and between
$\lambda_{F}^a(s_{k})$ and $\lambda_F^a(s_{k-1})$ are of order
$O\left(\frac{2k\pi^2}{4N^2}\right)=O\left(2k\epsilon^2\right)$
for $k=1,\dots, N.$ When the number of
atoms $N$ is sufficiently large, $\min_{0\le s\le 4}\lambda^a_{F}(s)$ can be used to
approximate $\min_{1\le k\le N}\lambda^a_{F}(s_{k})$
with ${1\le k\le
N}$ with
error at most of order
$O(\epsilon)$ since $N\eps =1.$

When $B_F<0$ and $N$ is sufficiently large, the minimum eigenvalue of $\delta^2\mathcal{E}^{a}(\mathbf{y}_{F})$ is no
longer $\lambda_{F}^{a}(s_1)$  and is given by
the following theorem.
\begin{theorem}\label{unstab}
Suppose that the hypotheses \eqref{EAMfunAssumption1} and $B_{F}< 0$
hold, and the number of atoms $N$ is sufficiently large. Then
$\lambda_{F}^{a}(s_{1})$ defined in Theorem~\ref{AtomStabThm} will no
longer be the minimum eigenvalue of the second variation
$\la\delta^{2}\mathcal{E}^{a}(\mathbf{y}_{F})\mathbf{u},\mathbf{u}\ra$.
Instead, the minimum eigenvalue will be given by $\lambda_{F}^a(s_{k^*})$
for some $s_{k^*}$, $1< k^*\le N$, that is either
equal to $4$ or close to
\[
s^*:=\frac{C_{F}-\sqrt{C_{F}^2-3B_{F}D_{F}}}{-3D_{F}}
\]
with difference of order
$O\left(2k^*\epsilon^2\right).$
\end{theorem}
\begin{proof}
Here we will briefly discuss the role of the coefficient $B_{F}$ and
leave the rigorous discussion of $\min_{0\le s\le 4}\lambda^{a}_{F}(s)$ under the condition $B_{F}<0$ to section~\ref{CompQNL1-2}.

The assumption $B_{F}\ge 0$ guarantees that
$u_{\ell}'=\sin(\epsilon \ell \pi)$ is the eigenfunction corresponding to
the minimum eigenvalue of $\delta^{2}\mathcal{E}^{a}(\mathbf{y}_{F})$
with respect to the norm $\|D\mathbf{u}\|_{\ell_{\epsilon}^{2}}.$
In fact, when $B_{F}<0$, we have ${\lambda^a_{F}}'(0)<0$ and thus
$\lambda^{a}_{F}(0)$ will be strictly larger than $\lambda^a_{F}(s^*)$.

We note that the condition $B_{F}\ge 0$ cannot be expected to generally
hold for EAM models when the nearest neighbor term $G''_{F}(\rho'_{F})^2>0$
dominates. We note, however, that generally $G'_{F}<0$ for $F<1$ ~\cite{Mishin}, in which case
$B_{F}\ge 0$ is more likely to hold for compressive
strains $F<1$.
\end{proof}
\begin{remark}\label{LargeNremark}
We would like to point out that when $N$ is small,
$\lambda_{F}^a(s_{1})$ may be still the minimum eigenvalue of
$\delta^2\mathcal{E}^a(\mathbf{y}_{F})$ even if $B_{F}<0$. This
is because $\lambda_F^a(s_{k})$ is defined on the discrete domain
$1\le k\le N$, so the continuous function $\lambda_{F}^a(s)$
is not a good approximation unless $N$ is
sufficiently large.
\end{remark}
\subsection{Stability of the Volume-Based and the Reconstruction-Based Local EAM Models.}
In this subsection, we will give stability estimations for the
volume-based and the reconstruction-based
local models, respectively.
\subsubsection{Stability of the Volume-Based Local EAM Model.}
We focus on the stability of the volume-based local model under a
uniform deformation $\mathbf{y}_{F}$. Using the equilibrium equation
\eqref{QCLEquilibrium1}, we obtain the second
variation $\delta^{2}\mathcal{E}^{c,v}(\mathbf{y}_{F})$
for any
$\mathbf{u}\in\mathcal{U}\setminus\{\mathbf {0}\}$
\[
\la\delta^{2}\mathcal{E}^{c,v}(\mathbf{y}_{F})\mathbf{u},\mathbf{u}\ra
=\left(\hat{A}_{F}+\tilde{A}_{F}\right)\|D\mathbf{u}\|^2_{\ell_{\epsilon}^2}
= A_{F}\|D\mathbf{u}\|^2_{\ell_{\epsilon}^2},
\]
where $\hat{A}_{F}$ and $\tilde{A}_{F}$ are defined in
\eqref{EAMCond1} and \eqref{PairAtomCond1}, respectively. It follows
that $\mathbf{y}_{F}$ is stable in the volume-based local model if and
only if $A_F:=\hat{A}_{F}+\tilde{A}_{F}>0$. We summarize this result in
the following theorem.
\begin{theorem}\label{QCLStabThm1}
Suppose that the hypotheses \eqref{EAMfunAssumption1} holds.
Then the uniform deformation $\mathbf{y}_{F}$ is stable in the volume-based local model \eqref{QCLenergy1}
if and only if $A_F:=\hat{A}_{F}+\tilde{A}_{F}>0$.
\end{theorem}
\begin{remark}\label{StabAtomvsQCL1}
Comparing the conclusions in Theorem \ref{AtomStabThm} and Theorem
\ref{QCLStabThm1}, we observe that when the hypothesis
\eqref{EAMfunAssumption2} is satisfied, the difference between the
minimum eigenvalues of the fully atomistic and the
volume-based local models is of order $O(\epsilon^2)$. This
result is the same as for the pair potential case
\cite{doblusort:qcf.stab}. However, when the assumptions fails, the
volume-based local model will be strictly more
stable than the fully atomistic model,
which will be discussed in the next remark.
\end{remark}
\begin{remark}\label{embed}
The assumption
\eqref{EAMfunAssumption2} is necessary for the validity of Theorem
\ref{AtomStabThm}. We now give an explicit example showing that the
uniform deformation can be strictly more stable for the volume-based
local model \eqref{QCLenergy1} than for the fully atomistic model
when \eqref{EAMfunAssumption2} fails.
We consider the case
\begin{equation}\label{EAMfunAssumption3}
\phi''_{2F}+G''_{F}\left(\rho'_{F}+2\rho'_{2F}\right)^2+G'_{F}2\rho''_{2F}>0.
\end{equation}
Then \eqref{EAMfunAssumption2} does not hold since it follows from
\eqref{EAMfunAssumption1} that
\begin{equation*}
\begin{split}
&\phi''_{2F}+G''_{F}\left[\left(\rho'_{F}\right)^2+20\left(\rho'_{2F}\right)^2
+12\rho'_{F}\rho'_{2F}\right]+G'_{F}2\rho''_{2F}\\
&\qquad=\left[\phi''_{2F}+G''_{F}\left(\rho'_{F}+2\rho'_{2F}\right)^2+G'_{F}2\rho''_{2F}\right]
+8G''_{F}\left(2\left(\rho'_{2F}\right)^2
+\rho'_{F}\rho'_{2F}\right)\\
&\qquad>0.
\end{split}
\end{equation*}

We define an oscillatory displacement $\tilde{\mathbf{u}},$ corresponding to the
$k=N$ eigenmode in the Fourier expansion \eqref{fourier}, by
\[
\tilde u_{\ell}:=(-1)^{\ell}\eps/(2\sqrt 2).
\]
Therefore,
\[
\tilde u_{\ell}'=(-1)^{\ell}/(\sqrt 2),\quad \|D\tilde{\mathbf{u}}\|_{\ell^2_{\epsilon}}=1,\quad
\tilde u_{\ell}''=(-1)^{\ell}(\sqrt 2)/\eps.
\]
From
\eqref{PairAtomStabEq1} and \eqref{AtomStabEq1} we can get
\begin{equation}\label{atomstab}
\begin{split}
\la\delta^2\mathcal{E}^a(\mathbf{y}_{F})\tilde{\mathbf{u}},\tilde{\mathbf{u}}\ra &=
\la\delta^2\tilde{\mathcal{E}}^a(\mathbf{y}_{F})\tilde{\mathbf{u}},\tilde{\mathbf{u}}\ra
+\la\delta^2\tilde{\mathcal{E}}^a(\mathbf{y}_{F})\tilde{\mathbf{u}},\tilde{\mathbf{u}}\ra\\
&=\epsilon\sum_{\ell=-N+1}^{N}G'_{F}2\rho''_{F}\half
+\left(\phi''_{F}+4\phi''_{2F}\right)\|D\tilde{\mathbf{u}}\|_{\ell^2_{\epsilon}}^2
+(-\epsilon^2\phi''_{2F})\|D^{(2)}\tilde{\mathbf{u}}\|^2_{\ell^2_{\epsilon}}\\
&=G'_{F}2\rho''_{F}+
\left(\phi''_{F}+4\phi''_{2F}\right)-4\phi''_{2F}=\phi''_{F}+G'_{F}2\rho''_{F}.
\end{split}
\end{equation}
Thus, we can obtain
\[
\inf_{{\mathbf{u}}\in\mathcal{U}\setminus\{\mathbf{0}\},\,\|D{\mathbf{u}}\|_{\ell^2_{\epsilon}}=1}
\la\delta^2\mathcal{E}^a(\mathbf{y}_{F}){\mathbf{u}},{\mathbf{u}}\ra\le
\phi''_{F}+G'_{F}2\rho''_{F}.
\]

On the other hand,
from Theorem~\ref{QCLStabThm1} we have that
 \begin{equation*}
\inf_{{\mathbf{u}}\in\mathcal{U}\setminus\{\mathbf{0}\},\,\|D{\mathbf{u}}\|_{\ell^2_{\epsilon}}=1}
\la\delta^2\mathcal{E}^{c,v}(\mathbf{y}_{F})\mathbf{u},\mathbf{u}\ra \equiv \tilde{A}_{F}+\tilde{A}_{F}
=4\left[\phi''_{2F}+G''_{F}\left(\rho'_{F}+2\rho'_{2F}\right)^2+G'_{F}2\rho''_{2F}\right]
+ \phi''_{F}+G'_{F}2\rho''_{F}.
\end{equation*}
Therefore, from \eqref{EAMfunAssumption3} and \eqref{atomstab} we have
\begin{align*}
\inf_{{\mathbf{u}}\in\mathcal{U}\setminus\{\mathbf{0}\},\,\|D{\mathbf{u}}\|_{\ell^2_{\epsilon}}=1}
\la\delta^2\mathcal{E}^{c,v}(\mathbf{y}_{F}){\mathbf{u}},{\mathbf{u}}\ra > \phi''_{F}+G'_{F}2\rho''_{F}
 \ge
\inf_{{\mathbf{u}}\in\mathcal{U}\setminus\{\mathbf{0}\},\,\|D{\mathbf{u}}\|_{\ell^2_{\epsilon}}=1}
\la\delta^2\mathcal{E}^a(\mathbf{y}_{F})\mathbf{u},\mathbf{u}\ra.
\end{align*}
This inequality indicates that when
the assumption \eqref{EAMfunAssumption2} fails, the uniform deformation
$\mathbf{y}_{F}$ can be unstable for the atomistic model, but still
stable for the volume-based local model.
\end{remark}
\subsubsection{Stability of the Reconstruction-Based Local EAM Model.}
In this case, we do a similar calculation for the reconstruction-based local model
and derive the second variation $\delta^{2}\mathcal{E}^{c,r}(\mathbf{y})$
from the equilibrium equation given by \eqref{QCLEquilibrium2}
\begin{align*}
\la\delta^{2}\mathcal{E}^{c,r} (\mathbf{y}_{F})\mathbf{u},\mathbf{u}\ra
=&\epsilon\sum_{\ell=-N+1}^{N}\left\{G''_{F}\left(\rho'_{F}+2\rho'_{2F}\right)^2\left(u'_{\ell}+u'_{\ell+1}\right)^2
+G'_{F}\left(\rho''_{F}+4\rho''_{2F}\right)\left[\left(u'_{\ell}\right)^2+\left(u'_{\ell+1}\right)^2\right]\right\}\notag\\
&\qquad+\epsilon\sum_{\ell=-N+1}^{N}\half
          \left\{\phi''_{F}\left[\left(u'_{\ell}\right)^2+\left(u'_{\ell+1}\right)^2\right]
+\phi''_{2F}\left[\left(4u'_{\ell}\right)^2+4\left(u'_{\ell+1}\right)^2\right]\right\}\\
=&\left[4G''_{F}\left(\rho'_{F}+2\rho'_{2F}\right)^2+2G'_{F}\left(\rho''_{F}+4\rho''_{2F}\right)
        +\phi''_{F}+4\phi''_{2F}\right]\|D\mathbf{u}\|^2_{\ell_{\epsilon}^{2}}\\
&\qquad\qquad
-\epsilon^2G''_{F}\left(\rho'_{F}+2\rho'_{2F}\right)^2\|D^{(2)}\mathbf{u}\|^2_{\ell_{\epsilon}^{2}}\\
=& A_{F}\|D\mathbf{u}\|^2_{\ell_{\epsilon}^{2}}+\epsilon^2\tilde{B}_{F}\|D^{(2)}\mathbf{u}\|^2_{\ell_{\epsilon}^{2}},
 \end{align*}
where $A_{F}$ is defined in \eqref{PairAtomCond1}
and $\tilde{B}_{F}$ is defined to be
\begin{equation}\label{QCL2Four}
\tilde{B}_{F}:=-G''_{F}\left(\rho'_{F}+2\rho'_{2F}\right)^2\le0.
\end{equation}
We recall that for the EAM-atomistic model, the coefficient $B_{F}$
of $\|D\mathbf{u}\|^2_{\ell_{\epsilon}^{2}}$ in \eqref{Four} is defined as
\[
B_{F}=-\left[\phi''_{2F}
+G''_{F}\left( \left(\rho'_{F}\right)^2+20\left(\rho'_{2F}\right)^2+12\rho'_{F}\rho'_{2F}\right)
+G'_{F}\left(2\rho''_{2F}\right)\right].
\]
Comparing $B_{F}$ with $\tilde{B}_{F}$ defined in \eqref{QCL2Four},
we find that
\[
B_{F}=\tilde{B}_{F}-\left[\phi''_{2F}
+G''_{F}\left( 16\left(\rho'_{2F}\right)^2+8\rho'_{F}\rho'_{2F}\right)
+G'_{F}\left(2\rho''_{2F}\right)\right].
\]
The assumption \eqref{EAMfunAssumption1} that $\phi''_{2F}<0$ implies
$B_{F}$ can be positive while $\tilde{B}_{F}$ is always
negative.

We similarly use the Fourier representation
\[
Du_{\ell}=\sum_{\substack {
k=-N+1\\ k\neq 0}}
^{N}\frac{c_{k}}{\sqrt{2}}\cdot\exp\left(i\,k\frac{\ell}{N}\pi\right)
\]
to analyze
the stability of $\delta^{2}\mathcal{E}^{c,r}(\mathbf{y}_{F}).$
Again, we exclude $k=0$ because of the mean zero condition of
$D\mathbf{u}$.
From the discrete orthogonality of the Fourier basis we have
\begin{equation}\label{QCL2StabFourierEq1}
\begin{split}
\la\delta^{2}\mathcal{E}^{c,r}(\mathbf{y}_{F})\mathbf{u},\mathbf{u}\ra &=
\sum_{\substack {
k=-N+1\\ k\neq 0}}
^{N} |c_{k}|^2\cdot \Bigg\{ A_F +\tilde{B}_F
\left[4\sin^{2}\left(\frac{k\pi}{2N}\right)\right]\Bigg\}.
\end{split}
\end{equation}
The eigenvalues
$\lambda_k^{c,r}$ of
$\la\delta^{2}\mathcal{E}^{c,r}(\mathbf{y}_{F})\mathbf{u},\mathbf{u}\ra$
with respect to the $\|D\mathbf{u}\|_{\ell_{\epsilon}^{2}}$ norm are
given by
\[
\lambda_k^{c,r}=\lambda_{F}^{c,r}(s_k)
\quad\text{for}\, k=1,\dots,N,
\]
where
\begin{equation*}
s_k=4\sin^2\left(\frac{k\pi}{2N}\right)\quad\text{and}\quad
\lambda_{F}^{c,r}(s):=A_{F}+\tilde{B}_{F}s.
\end{equation*}
The assumption \eqref{EAMfunAssumption1} implies that
$\tilde{B}_{F}\le 0$ always holds, so $\lambda_{F}^{c,r}(s)$ is
decreasing for $0\le s\le 4$ and the minimum eigenvalue of
$\delta^{2}\mathcal{E}^{c,r}(\mathbf{y}_{F})$ is achieved at $k=N$, i.e, $s_{N}=4$:
\begin{equation*}
\min_{1\le k\le
N}\lambda_{F}^{c,r}(s_{k})=\lambda_{F}^{c,r}(4)=A_{F}+4\tilde{B}_{F}=
2G'_{F}\left(\rho''_{F}+4\rho''_{2F}\right)  + \phi''_{F}
+4\phi''_{2F}.
\end{equation*}
The minimum eigenmode is given by the
oscillatory displacement
$\hat{u}'_{\ell}=-\hat{u}'_{\ell+1}$ since
\[
\la\delta^{2}\mathcal{E}^{c,r} (\mathbf{y}_{F})\hat{\mathbf{u}},\hat{\mathbf{u}}\ra
=\left[2G'_{F}\left(\rho''_{F}+4\rho''_{2F}\right)
        +\phi''_{F}+4\phi''_{2F}\right]\|D\hat{\mathbf{u}}\|^2_{\ell_{\epsilon}^{2}}.
\]
We thus have the following stability result for the reconstruction-based
local model.
\begin{theorem}\label{QCLStabThm2}
Suppose that the hypotheses \eqref{EAMfunAssumption1} holds.
Then the uniform deformation $\mathbf{y}_{F}$ is stable in the
reconstruction-based local model \eqref{QCLenergy2} if and
only if
\[2G'_{F}\left(\rho''_{F}+4\rho''_{2F}\right)  + \phi''_{F}
+4\phi''_{2F}>0.\]
\end{theorem}
\begin{remark}\label{QCL2unstable}
Comparing the conclusions in Theorem~\ref{AtomStabThm},
Theorem~\ref{QCLStabThm1}, and Theorem~\ref{QCLStabThm2}, we note
that when the assumption~\eqref{EAMfunAssumption2} holds, i.e.,
$B_{F}\ge 0$, the fully atomistic model is strictly more stable
than the reconstruction-based local model. The
difference between their minimum eigenvalues is $O(1)$,
not $O(\epsilon^2)$ as for the volume-based approximation.
When the assumption~\eqref{EAMfunAssumption2} fails, i.e. $B_{F}<0$,
the conclusion will be different, and we will rigorously analyze
this case in section~\ref{CompQNL1-2}.
\end{remark}
\section{Comparison of the Stability of the Atomistic and
Local EAM Models} \label{CompQNL1-2}
In this section, we would like
to give a full discussion of the sharp stability estimates for all of the
EAM models. Recall that the
eigenvalue function of $\delta^2\mathcal{E}^{a}(\mathbf{y}_{F})$ is
\[
\lambda^{a}_{F}(s_{k}):=A_{F}+B_{F}s_{k}+C_{F}s_{k}^2+D_{F}s_{k}^3
\quad \text{for}\quad
s_{k}=4\sin^2\left(\frac{k\pi}{2N}\right),\quad k=1,\dots,N,
\]
where the coefficients $A_{F}$, $B_{F}$, $C_{F}$ and $D_{F}$ are
given in the equation \eqref{AtomCoeff}.

To simplify the following analyses, the number of atoms $N$ is
assumed to be sufficiently large. Thus, we use the global minimum
of the continuous function $\lambda^{a}_{F}(s
):=A_{F}+B_{F}s +C_{F}s ^2+D_{F}s ^3$ for ${0\le s\le 4}$ to approximate
 $\min_{1\le k \le N}\lambda^{a}_{F}(s_{k})$. We note that their
difference is at most of order $O\left(2k\epsilon^2\right)\le
O(\epsilon)$.

We recall that
\begin{equation}\label{first}
\min_{0\le s \le 4}\lambda^{a}_{F}(s)=\lambda^{a}_{F}(0)\quad\text{if }B_F\ge 0.
\end{equation}
To find $\min_{0\le s\le 4}\lambda^a_{F}(s)$ when $B_F<0,$ we first evaluate $\lambda^a_{F}(s)$
at $s=0, 4$:
\begin{align*}
\lambda^a_{F}(0)=&A_{F}=4G''_{F}\left(\rho'_{F}+2\rho'_{2F}\right)^2+2G'_{F}\left(\rho''_{F}+4\rho''_{2F}\right)+\phi'_{F}+4\phi''_{2F},\\
\lambda^a_{F}(4)=&\phi''_{F}+2G'_{F}\rho''_{F}.
\end{align*}
We next compute the first and second derivatives of $\lambda^{a}_{F}(s)$, which are
\begin{align*}
{\lambda^a_{F}}'(s)=&B_{F}+2C_{F}s+3D_{F}s^2,\\
{\lambda^a_{F}}''(s)=&2C_{F}+6D_{F}s.
\end{align*}
Since ${\lambda^a_{F}}'(s)$ is a quadratic function,
we thus have two critical points of $\lambda^a_{F}(s)$ when the
coefficients satisfy
\[
C_{F}^2-3B_{F}D_{F}\ge 0
\quad\text{or equivalently}\quad\phi''_{2F}+2G'_{F}\rho''_{2F}\le \frac{1}{3}G''_{F}\left(\rho'_{F}-2\rho'_{2F}\right)^2.
\]
We can summarize the case when $B_F<0$ and $C_{F}^2-3B_{F}D_{F}\le0$ by
\begin{equation}\label{third}
\min_{0\le s \le 4}\lambda^{a}_{F}(s)=\lambda^{a}_{F}(4)\quad\text{if }B_F< 0\text{ and }
C_{F}^2-3B_{F}D_{F}\le 0.
\end{equation}

In the case $C_{F}^2-3B_{F}D_{F}> 0$ , the critical points are
\[
s_{1}=\frac{C_{F}-\sqrt{C^2_{F}-3B_{F}D_{F}}}{-3D_{F}} \quad
\text{and} \quad
s_{2}=\frac{C_{F}+\sqrt{C^2_{F}-3B_{F}D_{F}}}{-3D_{F}}.
\]
Since $D_F<0,$ $\lambda^a_{F}(s)$ will then have a local minimum
at $s^*=s_1$ and a local maximum at $s_2.$
The corresponding
local minimum value is
\begin{equation*}\label{AtomlocaMin}
\begin{split}
{\lambda^{a}_{F}}(s^*)= &A_{F}+B_{F}s^*+C_{F}(s^*)^2+D_{F}(s^*)^3\\
=& A_{F}+s^*\left(\frac{B_{F}}{2}-\frac{D_{F}}{2}\left(s^*\right)^2\right),
\end{split}
\end{equation*}
where we use ${\lambda^a_{F}}'(s^*)=0$ to get the last equality.
We can thus summarize all of the cases by
\begin{equation}\label{fourth}
\begin{split}
\min_{0\le s \le 4}\lambda^{a}_{F}(s)&=\lambda^{a}_{F}(0)\quad\text{if }B_F\ge 0,\\
\min_{0\le s \le 4}\lambda^{a}_{F}(s)&=\lambda^{a}_{F}(4)\quad\text{if }B_F< 0\text{ and }
C_{F}^2-3B_{F}D_{F}\le 0,\\
\min_{0\le s \le 4}\lambda^{a}_{F}(s)&=\min \{\lambda^{a}_{F}(s^*),\,\lambda^{a}_{F}(4)\}
 \quad\text{if }B_F< 0\text{ and }
C_{F}^2-3B_{F}D_{F}> 0.
\end{split}
\end{equation}

We note that the minimum eigenvalues of the volume-based and the
reconstruction-based local models are separately given by
the following expressions
\begin{equation}\label{compare}
\begin{split}
\min_{0\le s\le 4}\lambda^{c,v}_{F}(s)=&\lambda^{a}_{F}(0)=A_{F}
=4G''_{F}\left(\rho'_{F}+2\rho'_{2F}\right)^2
+2G'_{F}\left(\rho''_{F}+4\rho''_{2F}\right)+\phi''_{F}+4\phi''_{2F},\\
\min_{0\le s\le 4}\lambda^{c,r}_{F}(s)=&\lambda^{c,r}_{F}(4)=A_F-4G''_{F}\left(\rho'_{F}+2\rho'_{2F}\right)^2=
2G'_{F}\left(\rho''_{F}+4\rho''_{2F}\right)+\phi''_{F}
+4\phi''_{2F}.
\end{split}
\end{equation}
\subsection{The Volume-based Local EAM Model versus the Fully Atomistic EAM Model.}
We first compare the minimum eigenvalues of the volume-based local
and the fully atomistic models. Combining the
results of Theorem~\ref{AtomStabThm} and Theorem~\ref{unstab}, we
have the following theorem.
\begin{theorem}
The relation of the stability of the volume-based local model and
the fully atomistic model depends on the sign of
\[
B_{F}:= -\left[\phi''_{2F}+G''_{F}\Big((\rho'_{F})^2+20(\rho'_{2F})^2+12\rho'_{F}\rho'_{2F}\Big)+G'_{F}\left(2\rho''_{2F}\right)\right]
\]
and can be summarized as follows:
\begin{align*}
\min\limits_{0\le s\le 4}\lambda^{a}_{F}(s)=&\lambda^{a}_{F}(0)=\lambda^{c,v}_{F}\quad \text{if}\quad B_{F}\ge0,\\
\min_{0\le s\le 4}\lambda^{a}_{F}(s)=&\min\{\lambda^{a}_{F}(s^*),\lambda^{a}_{F}(4)\}
<\lambda^{a}_{F}(0)=\min_{0\le s\le 4}\lambda^{c,v}_{F} \quad \text{if}\quad B_{F}<0.
\end{align*}
\end{theorem}
This observation indicates that the set of stable uniform strains for the
volume-based local model {always includes} that for the fully atomistic EAM model.

\subsection{The Reconstruction-based Local EAM Model versus the Fully Atomistic EAM Model}
The relation of the minimum eigenvalues for $\delta^2\mathcal{E}^a(\mathbf{y}_{F})$
and $\delta^2\mathcal{E}^{c,r}(\mathbf{y}_{F})$
is more complicated.
We note that assumption \eqref{EAMfunAssumption1} implies
\begin{align*}
\lambda^a_{F}(0)=&A_{F}=4G''_{F}\left(\rho'_{F}+2\rho'_{2F}\right)^2
+2G'_{F}\left(\rho''_{F}+4\rho''_{2F}\right)+\phi''_{F}+4\phi''_{2F}\\
\ge& 2G'_{F}\left(\rho''_{F}+4\rho''_{2F}\right)+\phi''_{F}+4\phi''_{2F}=\min_{0\le s\le 4}\lambda^{c,r}_{F}(s),
\end{align*}
and we have
\[
\lambda^a_{F}(4)-\min_{0\le s\le 4}\lambda^{c,r}_{F}(s)=-4\left(\phi''_{2F}+G'_{F}\cdot2\rho''_{2F}\right).
\]
We thus conclude that if
\begin{equation}\label{EAMfunAssumption4}
\phi''_{2F}+G'_{F}\cdot2\rho''_{2F}\le 0,
\end{equation}
then
\[\lambda^a_{F}(4)\ge\min_{0\le s\le 4}\lambda^{c,r}_{F}(s).\]
 The equal sign
is achieved if and only if $\phi''_{2F}+G'_{F}\cdot2\rho''_{2F}=0$.
We also have the identity
\[
\phi''_{2F}+G'_{F}\cdot2\rho''_{2F}=-B_{F}-\left[G''_{F}\Big((\rho'_{F})^2+20(\rho'_{2F})^2
+12\rho'_{F}\rho'_{2F}\Big)\right].
\]

We next compare $\lambda^a_{F}(s^*)$ and $\min_{0\le s\le 4}\lambda^{c,r}_{F}(s)$.
The difference of these two is
\begin{align}%
\lambda^a_{F}(s^*)-&\min_{0\le s\le 4}\lambda^{c,r}_{F}(s)\nonumber\\
&= 4G''_{F}\left(\rho'_{F}+2\rho'_{2F}\right)^2+
s^*\left(\frac{B_{F}}{2}-\frac{D_{F}}{2}\left(s^*\right)^2\right)\nonumber\\
&= 4G''_{F}\left(\rho'_{F}+2\rho'_{2F}\right)^2\label{StabAtomvsQCL2Eq1}\\
&\quad+ \frac{C_{F}-\sqrt{C_{F}^2-3B_{F}D_{F}}}{-3D_{F}}\cdot
\frac{6B_{F}D_{F}-2C_{F}^2+C_{F}\left(C_{F}+\sqrt{C_{F}^2-3B_{F}D_{F}}\right)}{9D_{F}}\notag\\
&=4G''_{F}\left(\rho'_{F}+2\rho'_{2F}\right)^2+\frac{B_{F}C_{F}}{-9D_{F}}
 +\frac{2\left(C_{F}^2-3B_{F}D_{F}\right)\left(C_{F}-\sqrt{C_{F}^2-3B_{F}D_{F}}\right)}{27D_{F}^2}\notag\\
& \ge 4G''_{F}\left(\rho'_{F}+2\rho'_{2F}\right)^2+\frac{B_{F}C_{F}}{-9D_{F}}.\notag
\end{align}
According to the assumption \eqref{EAMfunAssumption4}, we
can use the
definition of $C_{F}$ and $D_{F}$ \eqref{AtomCoeff} to get
\begin{align*}
B_{F}=& -\left[\phi''_{2F}+G'_{F}\cdot 2\rho''_{2F}
+G''_{F}\left( (\rho'_{F})^2+12\rho'_{F}\rho'_{2F}+20(\rho'_{2F})^2 \right)\right]\\
\ge& -G''_{F}\left( (\rho'_{F})^2+12\rho'_{F}\rho'_{2F}+20(\rho'_{2F})^2 \right).
\end{align*}
We thus can obtain from the above inequality and the assumption \eqref{EAMfunAssumption1} that
\begin{align*}
\lambda^a_{F}(s^*)-&\min_{0\le s\le 4}\lambda^{c,r}_{F}(s)\\
&\ge 4G''_{F}\left(\rho'_{F}+2\rho'_{2F}\right)^2+\frac{B_{F}C_{F}}{-9D_{F}}\\
&\ge 4G''_{F}\left(\rho'_{F}+2\rho'_{2F}\right)^2
+\frac{G''_{F}\left(8\rho'_{2F}+2\rho'_{F}\right)
 \left( (\rho'_{F})^2+12\rho'_{F}\rho'_{2F}+20(\rho'_{2F})^2 \right)}{-9\rho'_{2F}}\\
&=2G''_{F}\frac{\left(\rho'_{F}+2\rho'_{2F}\right)\left(2\rho'_{2F}-\rho'_{F}\right)^2}{9\rho'_{2F}}\ge 0.
\end{align*}
Therefore, we have that
\begin{align}\label{StabAtomvsQCL2Summary1}
\begin{split}
\min_{0\le s\le 4}\lambda^a_{F}(s)=&\min\{\lambda^a_{F}(0),\lambda^a_{F}(s^*),\lambda^a_{F}(4)\}
= \min_{0\le s\le 4}\lambda^{c,r}_{F}(s)\quad\text{if}\quad \phi''_{2F}+G'_{F}\cdot2\rho''_{2F}=0,\\
\min_{0\le s\le 4}\lambda^a_{F}(s)=&\min\{\lambda^a_{F}(0),\lambda^a_{F}(s^*),\lambda^a_{F}(4)\}
> \min_{0\le s\le 4}\lambda^{c,r}_{F}(s)\quad\text{if}\quad\phi''_{2F}+G'_{F}\cdot2\rho''_{2F}<0.
\end{split}
\end{align}

Now let us turn to the case that the assumption \eqref{EAMfunAssumption4} fails, which means
\begin{equation}\label{EAMfunAssumption5}
\phi''_{2F}+G'_{F}\cdot2\rho''_{2F}> 0.
\end{equation}
In this case we have the opposite conclusion that the fully
atomistic model $\mathcal{E}^a(\mathbf{y})$ is strictly less stable
than the reconstruction-based local model
$\mathcal{E}^{c,r}(\mathbf{y})$.
From the condition
\eqref{EAMfunAssumption5}, we have
\begin{align*}
\lambda^a_{F}(4)-\min_{0\le s\le 4}\lambda^{c,r}_{F}(s)=&-4\left(\phi''_{2F}+G'_{F}\cdot2\rho''_{2F}\right)<0,\\
\text{i.e.,}\quad \lambda^a_{F}(4)<& \min_{0\le s\le 4}\lambda^{c,r}_{F}(s).
\end{align*}
Before comparing $\lambda^a_{F}(s^*)$ and $\min_{0\le s\le 4}\lambda^{c,r}_{F}(s)$,
we recall that $s^*$ exists if and only if
\[
\phi''_{2F}+2G'_{F}\rho''_{2F}\le \frac{1}{3}G''_{F}\left(\rho'_{F}-2\rho'_{2F}\right)^2.
\]
Thus, we actually consider the case that
\[
0<\phi''_{2F}+2G'_{F}\rho''_{2F}\le \frac{1}{3}G''_{F}\left(\rho'_{F}-2\rho'_{2F}\right)^2.
\]
We substitute
$B_{F}= -\left[\phi''_{2F}+G'_{F}\cdot 2\rho''_{2F}
+G''_{F}\left( (\rho'_{F})^2+12\rho'_{F}\rho'_{2F}+20(\rho'_{2F})^2 \right)\right]$
into the inequality \eqref{StabAtomvsQCL2Eq1} and get
\begin{align*}
\lambda^a_{F}(s^*)-&\min_{0\le s\le 4}\lambda^{c,r}_{F}(s)\\
 &\ge 4G''_{F}\left(\rho'_{F}+2\rho'_{2F}\right)^2
 +\frac{B_{F}C_{F}}{-9D_{F}}\notag\\
 &= 4G''_{F}\left(\rho'_{F}+2\rho'_{2F}\right)^2
 + \left(\frac{\phi''_{2F}+G'_{F}\cdot2\rho''_{2F}}{-9G''_{F}(\rho'_{2F})^2}
 +\frac{G''_{F}\left((\rho'_{F})^2+12\rho'_{F}\rho'_{2F}+20(\rho'_{2F})^2\right)}{-9G''_{F}(\rho'_{2F})^2}\right)C_{F}\\
 &= \frac{\left(\phi''_{2F}+G'_{F}\cdot2\rho''_{2F}\right)\left(8\rho'_{2F}+2\rho'_{F}\right)}{-9\rho'_{2F}}
 +2G''_{F}\frac{\left(\rho'_{F}+2\rho'_{2F}\right)\left(2\rho'_{2F}-\rho'_{F}\right)^2}{9\rho'_{2F}}.
\end{align*}
Since $\phi''_{2F}+2G'_{F}\rho''_{2F}\le \frac{1}{3}G''_{F}\left(\rho'_{F}-2\rho'_{2F}\right)^2$,
therefore
\begin{align*}
\lambda^a_{F}(s^*)-&\min_{0\le s\le 4}\lambda^{c,r}_{F}(s)\\
&\ge \frac{\left(\phi''_{2F}+G'_{F}\cdot2\rho''_{2F}\right)\left(8\rho'_{2F}+2\rho'_{F}\right)}{-9\rho'_{2F}}
 +2G''_{F}\frac{\left(\rho'_{F}+2\rho'_{2F}\right)\left(2\rho'_{2F}-\rho'_{F}\right)^2}{9\rho'_{2F}}\\
& \ge  \frac{\frac{1}{3}G''_{F}\left(\rho'_{F}-2\rho'_{2F}\right)^2\left(8\rho'_{2F}+2\rho'_{F}\right)}{-9\rho'_{2F}}
 +2G''_{F}\frac{\left(\rho'_{F}+2\rho'_{2F}\right)\left(2\rho'_{2F}-\rho'_{F}\right)^2}{9\rho'_{2F}}\\
& = \frac{1}{3}G''_{F}\left(\rho'_{F}-2\rho'_{2F}\right)^2\frac{4\rho'_{2F}+4\rho'_{F}}{9\rho'_{2F}}\ge 0.
\end{align*}
Hence, when $s^*$ exists and the assumption \eqref{EAMfunAssumption5} holds, we have
\[
\min_{0\le s\le 4}\lambda^a_{F}(s)=\min\{\lambda^a_{F}(0),\lambda^a_{F}(s^*),\lambda^a_{F}(4)\}
=\lambda^a_{F}(4)<\min_{0\le s\le 4}\lambda^{c,r}_{F}(s).
\]
When $s^*$ does not exist,
we can immediately get that
\[
\min_{0\le s\le 4}\lambda^a_{F}(s)=\min\{\lambda^a_{F}(0),
\lambda^a_{F}(4)\}=\lambda^a_{F}(4)<\min_{0\le s\le 4}\lambda^{c,r}_{F}(s).
\]
We now combine this result with \eqref{StabAtomvsQCL2Summary1} and
summarize the stability relation between the fully atomistic model
and the reconstruction-based local model by the following
theorem.
\begin{theorem}
The relation between the stability of the reconstruction-based local
 model  and the atomistic model
depends on the sign of  $\phi''_{2F}+G'_{F}\cdot2\rho''_{2F}$ and is given by
\begin{align}\label{StabAtomvsQCL2Summary2}
\begin{split}
\min_{0\le s\le 4}\lambda^{a}_{F}(s)
=&\lambda^{a}_{F}(4)<\min_{0\le s\le 4}\lambda^{c,r}_{F}(s)\quad\text{if}\quad\phi''_{2F}+G'_{F}\cdot2\rho''_{2F}>0,\\
\min_{0\le s\le 4}\lambda^a_{F}(s)=&\min\{\lambda^a_{F}(0),\lambda^a_{F}(s^*),\lambda^a_{F}(4)\}
= \min_{0\le s\le 4}\lambda^{c,r}_{F}(s)\quad\text{if}\quad \phi''_{2F}+G'_{F}\cdot2\rho''_{2F}=0,\\
\min_{0\le s\le 4}\lambda^a_{F}(s)=&\min\{\lambda^a_{F}(0),\lambda^a_{F}(s^*),\lambda^a_{F}(4)\}
> \min_{0\le s\le 4}\lambda^{c,r}_{F}(s)\quad\text{if}\quad\phi''_{2F}+G'_{F}\cdot2\rho''_{2F}<0.\\
\end{split}
\end{align}
\end{theorem}
We note from the theorem that the reconstruction-based local model {can be less stable than} the fully atomistic model,
which might cause stability problems when constructing a coupling method.

\section{Conclusion.}
{In this paper, we give precise estimates for the lattice stability
of atomistic chains modeled by the fully
atomistic EAM model and the volume-based and the
reconstruction-based local approximations. We identify
 the critical assumptions for the pair potential, the electron density
function, and the embedding function to study lattice stability.
We find that both the volume-based local model
and the reconstruction-based local model can give O$(1)$ errors for the
critical strain.  The critical strain predicted by the volume-based model is
always larger than that predicted by the atomistic model, but the critical
strain for reconstruction-based models can be either larger or smaller than
that predicted by the atomistic model.

Further research is needed to determine the significance of these results
for multidimensional lattice stability and for atomistic-to-continuum coupling
methods that couple an atomistic region with a volume-based local region
through a reconstruction-based local region~\cite{Shimokawa:2004,E:2004}.
}


\end{document}